\documentclass[12pt,english]{article}
\usepackage{geometry}
\usepackage{color}
\usepackage{amsmath}
\usepackage{amssymb}
\usepackage{amsthm}

\def\E{\hskip.15ex\mathrm{E}\hskip.10ex}
\def\P{\mathrm{P}}

\def\phi{\varphi}

\newif\ifabc
\abcfalse



\usepackage{babel}

\newtheorem{theorem}{Theorem}
\newtheorem{lemma}{Lemma}
\newtheorem{remark}{Remark}
\newtheorem{definition}
{Definition}
{Exercise}
{Assumption}
\newtheorem{corollary}
{Corollary}
\newtheorem{proposition}
{Proposition}

\newtheorem{example}
{Example}

\begin{document}
\global\long\def\E{\mathbb{E}}
\global\long\def\P{\mathbb{P}}
\global\long\def\N{\mathbb{N}}
\global\long\def\ind{\mathbb{I}}

\title{
{\normalsize\tt\hfill\jobname.tex}\\
On convergence rate for homogeneous Markov chains\\
}
\author{A.Yu. Veretennikov\footnote{ University of Leeds, UK; National Research University Higher School of Economics, and Institute for Information Transmission Problems, Moscow, Russia, email: a.veretennikov @ leeds.ac.uk. Supported by 
the grant RSF 17-11-01098}, 
M.A. Veretennikova\footnote{National Research University Higher School of Economics, Moscow, Russia, email: mveretennikova @ hse.ru. Supported by the  grant RSF 17-11-01098}
}

\maketitle
\begin{abstract}
This paper 
studies 
improved rates of convergence for ergodic homogeneous Markov chains. In comparison to the earlier papers the setting is also generalised to the case without a unique dominated measure.
Examples are provided where the new bound is compared with 
the classical Markov -- Dobrushin inequality and with the second eigenvalue of the transition matrix for finite state spaces. 
\end{abstract}


\section{Introduction} 
It is well-known that for a discrete irreducible acyclic (also called primitive) homogeneous Markov chain $(X_n, n\ge 0)$ with a finite state space $\cal S$ there is a unique stationary distribution $\mu$, and the distribution of the chain $\mu_n={\cal L}(X_n)$ admits an exponential bound uniformly with respect to the initial distribution $\mu_0$:
\begin{equation}\label{finiteMD}
\|\mu_n - \mu\|_{TV} \le 2 (1-\kappa)^n, \quad \forall n.
\end{equation}
This bound converges exponentially fast to zero if  
\begin{equation*}
\kappa = \min_{i,i'}\sum_{j\in {\cal S}} p_{ij}\wedge p_{i'j} >0, 
\end{equation*}
where $p_{ij}$ are the transition probabilities, see, e.g., \cite{Gnedenko}.  Here $\|\cdot\|_{TV}$ is the distance of total variation between measures. A similar bound exists for general Markov chains, too (cf., among many other sources, \cite{Veretennikov17}). Another well-known method of exponential estimates is related to eigenvalues of the transition probability matrix, see \cite[Chapter XIII, (96)]{Gantmacher}, \cite[Theorem 1.2]{Seneta1}. However, in the general situation it suits special cases where the eigenvectors (eigenfunctions) of the transition operator form a full basis in $L_2({\cal S})$. This approach works especially well  in the simplest case of $|{\cal S}| = 2$ (see, e.g., the Proposition \ref{pros2} below) and for symmetric cases which correspond to reversible Markov chains (see, e.g., \cite{Stroock-Diaconis}). 
In \cite{BV} and \cite{Veretennikov17} it was prompted that a certain markovian coupling (see \cite{Griffeath, Lindvall, Thorisson, Vaserstein} et al.) construction may in some cases provide a better convergence rate than that of (\ref{finiteMD}); note that this approach is applicable to non-reversible cases, too. Yet, neither paper provided examples \cite{BV,Veretennikov17}. Another small issue in both \cite{BV} and \cite{Veretennikov17} was a reasonable but still slightly restrictive an assumption of existence of a unique reference measure (in \cite{BV} the  Lebesgue one) for all transition kernels $Q(x, \cdot)= P(X_1 \in \cdot | X_0)|_{X_0=x}$. This assumption is dropped in the present paper. 

We emphasize that the eigenvalue method provides the best asymptotic convergence rate for finite transition matrices. However, the advantage of bound (\ref{finiteMD}) is that it is valid in general state spaces, see any texbook treating such a case, or, for example, \cite{Veretennikov17}, or chapter \ref{sec:erg_gen} below. Beside this generality, here some -- hopefully convincing -- examples are provided, which are inevitably finite (so far). 

These examples demonstrate that the new bound is effective -- unlike, e.g., the remarkable  Doeblin--Doob one \cite[Theorem V.5.6]{Doob53}. The only drawback of the latter is that the ``condition D'' from \cite{Doob53}, which is nowadays called the Doeblin-Doob condition, does not imply any computable bound on the constants in the resulting convergence rate, so this rate may actually be arbitrarily slow, albeit being exponential. The new approach  provides the same or a better bound in comparison to (\ref{finiteMD}).   
The authors are working on further examples for infinite state spaces.

The rest of the paper consists of two sections: 2 -- the presentation of the method from \cite{BV} and \cite{Veretennikov17} in general state spaces and under less restrictive assumptions (without a unique reference measure); 3 -- examples of both sorts, i.e., where the new bound is better and where it gives basically the same result as the classical MD inequality. 

\section{MD Ergodic Theorem and coupling, general case}\label{sec:erg_gen} 
We consider a homogeneous Markov process (MP) in discrete time $(X_n, \, n\ge 0)$ on a general state space $S$ with a topology and with its Borel sigma-algebra ${\cal B}(S)$.  
%
The following notations from the theory of Markov processes will be accepted (cf. \cite{Dynkin}): the index $x$ in $\mathbb E_x$ or $\mathbb P_x$ signifies the expectation or the probability measure related to the non-random initial state of the process $X_0$. This initial state may also be  random with some distribution $\mu$, in which case notations  $\mathbb E_\mu$ or $\mathbb P_\mu$ can be used. 

If  the state space $S$ is finite, then $|S|$ denotes the number of its elements and $\cal P$ stands for the transition matrix $\left(p_{ij}\right)_{1\le i,j \le |S|}$ of the process. Concerning the history of the MD bound see  \cite{Doob53, Dynkin, EthierKurtz, Seneta1}, et al.

For the most well-known inequality (\ref{exp_bd3}) of the Proposition \ref{thm_erg2} about the classical convergence rate bounds given below for the reader's convenience, it is assumed that the {\em Markov--Dobrushin constant} is positive:
\begin{equation}\label{MD}
 \kappa := \inf_{x,x'} \int \left(\frac{P_{x'}(1,dy)}{P_x(1,dy)}\wedge 1 \right)P_x(1,dy) > 0. 
\end{equation}
In the main result of this paper -- Theorem \ref{lastthm} in  the section -- this condition will be dropped. 
Note that here $\displaystyle \frac{P_{x'}(1,dy)}{P_x(1,dy)}$ is understood in the sense of the density of the absolute continuous components. For brevity we will be using a simplified notation $P_x(dz)$ for $P_x(1,dz)$. Note that integration of the kernel $P_x(1,dz)$ with any Borel function  is Borel measurable with respect to $x$ which is a standard requirement in Markov processes \cite{Dynkin}. An analogous measurability with respect to the pair $(x,x')$ will also be valid for an integral of a Borel function with the measure $\Lambda_{x,x'}$ defined by the formula
\begin{equation*}
\Lambda_{x,x'}(dz) := P_x(1,dz) + P_{x'}(1,dz), 
\end{equation*}  
due to linearity. 
Note that $\Lambda_{x,x'}(dz) = \Lambda_{x',x}(dz)$.
\begin{lemma}\label{newMD}
The following representation for the constant from (\ref{MD}) holds true: 
\begin{equation}\label{betterMD}
 \kappa = \inf_{x,x'} \int \left(\frac{P_{x'}(1,dy)}{\Lambda_{x,x'}(dy)}\wedge \frac{P_{x}(1,dy)}{\Lambda_{x,x'}(dy)} \right)\Lambda_{x,x'} (dy) \;\; (>0). 
\end{equation}
\end{lemma}

\noindent
{\em Proof } -- straightforward.  

~

Here is the key notion in the following presentation: denote 
\[
\kappa(x,x') : = \int \left(\frac{P_{x'}(1,dy)}{P_x(1,dy)}\wedge 1 \right)P_x(1,dy).
\]

\begin{lemma}\label{jj}
For any $x,x'\in S$,
\[
\kappa(x,x') = \kappa(x',x).
\]
\end{lemma}
\noindent
{\em Proof.} 
We have, 
 
\begin{eqnarray}
& \displaystyle \kappa(x',x) = \int \left(\frac{P_{x'}(1,dy)}{P_x(1,dy)}\wedge 1 \right)P_x(1,dy) 
= \int \left(\frac{P_{x'}(1,dy)}{P_x(1,dy)}\wedge 1 \right)\frac{P_x(1,dy)}{\Lambda_{x,x'}(dy)} \Lambda_{x,x'}(dy)
 \nonumber \\ \nonumber \\
& \displaystyle = \int \left(\frac{P_{x'}(1,dy)}{\Lambda_{x,x'}(dy)} \wedge \frac{P_x(1,dy)}{\Lambda_{x,x'}(dy)} \right) \Lambda_{x,x'}(dy). \label{newgamma}
\end{eqnarray}
The  Lemma \ref{jj} follows. 

\begin{definition}
We call Markov processes satisfying the condition (\ref{MD}) {\em Markov--Dobrushin's or {\bf MD}-processes}.
\end{definition}

In the scenario of finite homogeneous chains this condition (\ref{MD}) was introduced by Markov himself \cite{Markov}; later on, for non-homogeneous Markov processes its analogue was used by Dobrushin \cite{Dobrushin}. It was suggested by Seneta 
to call it the Markov--Dobrushin condition. Note that in all cases $\kappa \le 1$. The case $\kappa = 1$ corresponds to the i.i.d. sequence $(X_n, n\ge 0)$. In the opposite extreme situation where the transition kernels are singular for different $x$ and $x'$, we have $\kappa = 0$. 
The MD-condition (\ref{MD}) -- as well as (\ref{betterMD}) -- is most useful because it provides an {\bf effective} quantitative upper bound for the convergence rate of a Markov chain towards its (unique) invariant measure in total variation metric. The following theorem is classical: the bound (\ref{exp_bd3}) can be found in most textbooks on ergodic Markov chains. 

\begin{proposition}\label{thm_erg2}
Let the assumption (\ref{betterMD}) 
hold true. Then the process $(X_n)$ is ergodic, i.e., there exists  a limiting probability measure $\mu$, which is stationary and such that the uniform bound is satisfied for every $n$, 
\begin{equation}\label{exp_bd3}
 \sup_{x}\sup_{A\in {\cal B}(S)} |P_x(n,A) - \mu(A)| \le (1-\kappa)^{n}.
\end{equation}
\end{proposition}
Equivalently, for any bounded Borel measurable function $g$ on $S$, 
\begin{equation}\label{exp_bd6}
\sup_x |\langle g,\mu^x_n\rangle - \langle g, \mu \rangle | \le |g|_B (1-\kappa^{})^{n},
\end{equation}
where $|g|_B = \sup_x |g(x)|$.

~

The following folklore lemma answers the following question: suppose we have two distributions, which are not singular and their ``common area'' equals some positive constant $\kappa$. Is it possible to realise these two distributions on one  probability space so that the two corresponding random variables {\em coincide} with probability equal exactly $\kappa$? 

\begin{lemma}[``Of two random variables'']\label{odvuh}
Let $\xi^{1}$ and $\xi^2$ be two random variables on their (different, without loss of generality) probability spaces $(\Omega^1, {\cal F}^1, \mathbb P^1)$ and $(\Omega^2, {\cal F}^2, \mathbb P^2)$ and with densities $p^1$ and $p^2$ with respect to some reference measure $\Lambda$, correspondingly.  Then, if 
\begin{equation*}
q := \int \left(p^1(x)\wedge p^2(x)\right) \Lambda(dx) > 0, 
\end{equation*}
then there exists one more probability space $(\Omega, {\cal F}, \mathbb P)$ and two random variables on it, $\eta^1, \eta^2$, such that 
\begin{equation*}
{\cal L}(\eta^j) ={\cal L}(\xi^j), \; j=1,2, \quad \& \quad \mathbb  P(\eta^1 = \eta^2) = q. 
\end{equation*}
\end{lemma}
For the proof of this simple well-known fact see, e.g., \cite{Veretennikov17}. Here $q$ plays the role of $\kappa$.

In this section it is explained how to apply the general  coupling method  to Markov chains in general state spaces $(S, {\cal B }(S))$. Various presentations of this method may be found in \cite{Kalash, Lindvall, Nummelin, Thorisson, Vaserstein},  et al. This section follows the lines from \cite{BV}, which, in turn, is based on \cite{Vaserstein}. Note that in \cite{BV} the state space was $\mathbb R^1$; however, in $\mathbb R^d$ all the formulae remain the same. 

Now, let us consider two versions of the same Markov process $(X^1_n), (X^2_n)$  with two initial distributions $\mu_0^1$ and  $\mu_0^2$ respectively (this does not exclude the case of non-random initial states). Denote 
\begin{equation*}
\kappa_0 = \kappa(\mu_0^1, \mu_0^2):=
\int \left(\frac{\mu_0^1(dy)}{\mu_0^2(dy)}\wedge 1 \right)\mu_0^2(dy).
\end{equation*}
It is clear that $0\le \kappa_0\le1$ similarly to $\kappa(u,v)$ for all $u,v$.  
If $X^1_0$ and $X^2_0$ have different distributions, then obviously $\kappa_0<1$. Otherwise, we  have
$X^1_n\stackrel{d}{=}X^2_n$ (equality in distribution) for all $n$, and the coupling can be made trivially by letting  $\widetilde X^1_n= \widetilde
X^2_n:=X^1_n$.

Let us introduce a new, vector-valued {\bf Markov process} $\left(\eta^1_n,\eta^2_n,\xi_n,\zeta_n\right)$.  If $\kappa_0=0$ then we set
\begin{equation*}
\eta^1_0:=X^1_0,\; \eta^2_0:=X^2_0,\; \xi_0:=0,\; \zeta_0:=1.
\end{equation*}
Otherwise, if $0<\kappa(0)<1$, then we apply the Lemma \ref{odvuh} with $\Lambda = \Lambda_{x,x'}$ conditional on $(X^1_0, X^2_0)=(x,x')$ to the random variables $X^1_0$ and $X^2_0$ so as to create the random variables 
$\eta^1_0$, $\eta^2_0$, $\xi_0$ and $\zeta_0$ (they correspond to $\eta^1, \eta^2, \xi$, and $\zeta$ in the Lemma); formally, this can be done if $\kappa_0=1$, too. Further, assuming that the random variables $\left(\eta^1_n,\eta^2_n,\xi_n,\zeta_n\right)$ have been determined for some $n$, let us show how to construct the vector $\left(\eta^1_{n+1},\eta^2_{n+1},\xi_{n+1},\zeta_{n+1}\right)$. For this aim, we define the transition probability density $\phi$ with respect to the measure $\Lambda_{x^1, x^2} \times \Lambda_{x^1, x^2} \times \Lambda_{x^1, x^2}\times (\delta_0 + \delta_1)$ for this vector-valued process as follows,
\begin{equation}\label{process_eta}
\phi(x,y):=\phi_1(x,y^1)\phi_2(x,y^2)\phi_3(x,y^3) \phi_4(x,y^4),
\end{equation}
where $x=(x^1,x^2,x^3,x^4)$, $y=(y^1,y^2,y^3,y^4)$, and if
 $0<\kappa(x^1,x^2)<1$, then
\begin{align}
&\displaystyle \phi_1(x,u):=\frac{p(x^1,u)-p(x^1,
u)\wedge p(x^2,u)}{1-\kappa(x^1,x^2)}, \quad
\phi_2(x,u):=\frac{p(x^2,u)-p(x^1,u)\wedge p(x^2,u)}{1-\kappa(x^1,x^2)},\label{phi_12}
 \\\nonumber\\
&\displaystyle \phi_3(x,u):=1(x^4=1)\frac{p(x^1,u)\wedge
 p(x^2,u)}{\kappa(x^1,x^2)}+1(x^4=0)p(x^3,u),\label{phi_3}
 \\\nonumber\\ 
&\displaystyle \phi_4(x,u):=1(x^4=1)\left(\delta_1(u)(1-\kappa(x^1,x^2))+ 
\delta_0(u)\kappa(x^1,x^2)\right) +1(x^4=0)\delta_0(u)\label{phi_4}, 
\end{align}
where $\delta_i(u)$ is the Kronecker symbol, $\delta_i(u) = 1(u=i)$, or, in other words, the delta measure concentrated at state $i$. The case $x^4=0$ signifies coupling which has already been realised at the previous step, and $u=0$ means successful coupling at the transition.  
In the degenerate cases, if $\kappa(x^1,x^2)=0$ (coupling at the transition is impossible), then instead of (\ref{phi_3}) we set,  e.g., 
$$
\phi_3(x,u):=1(x^4=1)p(x^3,u) + 1(x^4=0)p(x^3,u) = p(x^3,u), \qquad 
$$ 
and if $\kappa(x^1,x^2)=1$, then instead of (\ref{phi_12}) we may set 
$$
\phi_1(x,u)=\phi_2(x,u):= p(x^1,u). \qquad 
$$ 
The formula (\ref{phi_4}) which defines \(\phi_4(x,u)\) can be accepted in all cases.

\begin{lemma}[\cite{BV}]\label{lemma:2}Let $X^1_0=x, X^2_0=x'$, 
and the random variables $\widetilde X^1_n$ and $\widetilde X^2_n$, for $n\in\mathbb{Z}_+$ be defined   by the following formulae:
\begin{align*}
\widetilde X^1_n:=\eta^1_n 1(\zeta_n=1)+\xi_n 1(\zeta_n=0), \quad 
\widetilde X^2_n:=\eta^2_n 1(\zeta_n=1)+\xi_n 1(\zeta_n=0).
\end{align*}
Then $\kappa_0 = \delta_{x,x'}$, and
\[
\widetilde X^1_n\stackrel{d}{=}X^1_n, \;\;\widetilde
 X^2_n\stackrel{d}{=}X^2_n, \quad \mbox{for all $n\ge 0$.}
\]
Moreover, 
\[
\widetilde X^1_n=\widetilde X^2_n, \quad \forall \; n\ge
n_0(\omega):=\inf\{k\ge0: \zeta_k=0\}, 
\] 
and
\begin{equation}\label{estimate}
\P_{x,x'}(\widetilde X^1_n\neq \widetilde
 X^2_n)\le(1-\kappa_0)\, \E_{x,x'}\prod_{i=0}^{n-1}
 (1-\kappa(\eta^1_i,\eta^2_i)).
\end{equation}
Moreover, $\left(\widetilde X^1_n\right)_{n\ge 0}$ and $\left(\widetilde X^2_n\right)_{n\ge 0}$ are both  homogeneous Markov processes, and 
\begin{equation*}
\left(\widetilde X^1_n\right)_{n\ge 0}\stackrel{d}{=}\left(X^1_n\right)_{n\ge 0}, \quad \& \quad 
\left(\widetilde X^2_n\right)_{n\ge 0}\stackrel{d}{=}\left(X^2_n\right)_{n\ge 0}.
\end{equation*}
\end{lemma}

\noindent
In \cite{BV} it was additionally assumed   in the Lemma \ref{lemma:2} that  there is a unique dominating measure for all transition kernels. This  assumption is not necessary and the proof runs without other  changes with $\Lambda_{x,x'}$ at each step instead of a unique dominating measure $\Lambda$. 

\begin{remark}\label{rem0}
If the initial values $X^1_0, X^2_0$ are distributed with laws $\mu^1_0, \mu_0^2$, correspondingly, then the following corollary  of the inequality (\ref{estimate}) holds, 
\begin{equation}\label{estimate2}
\P_{\mu^1_0,\mu^2_0}(\widetilde X^1_n\neq \widetilde
 X^2_n)\le(1-\kappa_0(\mu^1_0,\mu^2_0))\, \E_{\nu^1_0,\nu^2_0}\prod_{i=0}^{n-1}
 (1-\kappa(\eta^1_i,\eta^2_i)),
\end{equation}
where $\nu^j_0$ is the distribution of $\eta^j_0$, $j=1,2$. 

\end{remark}

~

Denote $\mbox{diag}(S^2): =  \{x\in S \times S: \, x^1=x^2\}$ and $\hat S^2 := S \times S \setminus \mbox{diag}(S^2)$. Consider the process $X_n := (\tilde X^1_n, \tilde X^2_n)$ on $S^2$ and its conditional version $\hat X_n := (\tilde X^1_n, \tilde X^2_n)|\tilde X^1_n \not = \tilde X^2_n$ on $\hat S^2$.
Note that $\hat X_n$ is also a Markov process. Indeed, its probabilities of states at time $n+1$ only depend on the state at time $n$, or, more precisely, given this state they do not depend on the past before $n$: in comparison to $X_n$ they are just recomputed so as to satisfy the normalisation condition.

~

Now using the function $\kappa(x)$, let us introduce an  operator $V$ acting on a (bounded, Borel measurable) function $h$ on the space $\hat S^2$ as follows: for $x=(x^1, x^2)\in  \hat S^2$,
\begin{equation}\label{V}
Vh(x) := (1-\kappa(x^1,x^2)) \mathbb E_{x^1,x^2}h(\hat X_1) \equiv \exp(\psi(x))\mathbb E_{x^1,x^2}h(\hat X_1),  
\end{equation}
where in the last expression $\psi(x):= \ln (1-\kappa(x^1,x^2))$ (assume $\ln 0 = -\infty$). From the estimate (\ref{estimate2}) it follows, 
\begin{equation}\label{l2}
\P_{\mu^1_0,\mu^2_0}(\widetilde X^1_n\neq \widetilde X^2_n)\le 
(1-\kappa_0(\mu^1_0,\mu^2_0)) \int V^n 1(x)\nu^1_0(dx^1)\nu^2_0(dx^2). 
\end{equation}
Note that by definition (\ref{V}), for the non-negative operator $V$ (which transforms any non-negative function into a non-negative one) its sup-norm $\|V\| = \|V\|_{B,B}:=\sup\limits_{|h|_B\le 1} |Vh|_B $ equals $\sup\limits_{x} V1(x)$, where 
$1=1(x)$ is the function on $\hat S^2$ identically equal to one. In our case 
$$
r(V) \le \|V\| = \sup\limits_{x} V1(x) = \sup_x \exp(\psi(x))\equiv \sup_x (1-\kappa(x)) = 1-\kappa, 
$$ 
due to the well-known inequality (see, for example, \cite[\S 8]{KLS}). 

Further, if the operator $V$ was compact (cf. \cite{KLS}), e.g., in the space $C(S)$ -- as is always true for finite $S$ -- then from the generalisation of the Perron--Frobenius Theorem (see, for example, \cite[\S 9, Theorem 9.2]{KLS}) it would follow (see, e.g., \cite[(7.4.10)]{FW}) that 
\begin{equation*}
\lim_{n\to\infty} \frac1n \, \ln V^n 1(x)  = \ln r(V). 
\end{equation*}
We also have, 
\(\sup_x V^n 1 (x) = \|V^n\|\) and  \(0\le  V^n 1 (x) \le  \|V^n\|\), for any \(x\in \hat S^2\). 
So from the Gelfand formula, 
\begin{equation}\label{l11}
\limsup_{n\to\infty} (V^n 1 (x))^{1/n} \le \lim_{n\to\infty} \|V^n\|^{1/n} =  r(V) \le \|V\|. 
\end{equation}
The latter inequality holds without any assumptions of compactness on $V$.
The assertions (\ref{l2}) and (\ref{l11}) together lead to the following result. 
\begin{theorem}\label{lastthm}
In all cases, 
\begin{equation*}
\limsup\limits_{n\to\infty} \frac1n \ln \| P_{x^1}(n,\cdot) - \mu(\cdot)\|_{TV} 
\le \limsup_{n\to\infty} \frac1n \, \ln \int V^n 1(x^1, x^2)\mu(dx^2)
\le \ln r(V).  
\end{equation*}
\end{theorem}
\begin{remark}\label{rem1}
If the analogue of the ergodicity condition (\ref{MD}) {\em is valid for the process $\hat X_n$,} 
then the strict inequality
\[
r(V) < \|V\| = 1-\kappa
\] 
holds iff $1-\kappa(\cdot)$ is {\em not a constant $\hat\mu$-a.s.,} where $\hat\mu$ is the unique invariant measure for $\hat X_n$.

\end{remark}


\section{Examples: Finite $S$}
\begin{proposition}\label{pros2}
Let us consider a $2 \times 2$ matrix $\cal{P}$ with elements $p_{1,1}=a, p_{1,2}=1-a, p_{2,1}=(1-b), p_{2,2}=b$, assuming $0<a,b<1$. 	Then $r=1-\kappa = |\lambda_{2}|$, where $\lambda_{2}$ is the smaller eigenvalue of the  matrix $\cal{P}$. 
\end{proposition}

{\em Proof} -- is straightforward and will be shown in the full version of the paper.  

~

In the examples below we compare the asymptotics provided by three methods: Markov--Dobrushin's, the one based on the second eigenvalue, and the one proposed in the Theorem \ref{lastthm} above. Although, in all examples the eigenvalues can be computed by hand, a code was written which can be applied to any Markov chain in a finite state; this code will be shown in the full version of the paper along with further multidimensional examples. 
\begin{example}[the three methods provide  similar rate bounds]
Consider the MP with the state space $S = \{1,2\}$ and a transition matrix
\[
{\cal P} = \left(
\begin{array}{l l}
0.65 & 0.35 \\
0.35 & 0.65
\end{array}
\right).
\]
We have, $\hat S^2 = \{(1,2), (2,1)\}$, and 
\[
V = (1-0.7) \times 
\left(
\begin{array}{l l}
0.65 - 0.35 & 0 \\
0 & 0.65 - 0.35
\end{array}
\right)\times \frac{1}{0.3} = 0.3 \times 
\left(
\begin{array}{l l}
1 & 0 \\
0 & 1
\end{array}
\right)\,.
\]
So we compute, 
\[
1-\kappa = 0.3   =  r(V), 
\]
as expected. We have, 
\[
\lambda_2 = 0.3 = r(V), 
\]
again as expected since here $\lambda_2 = 1-\kappa$.  
This Markov chain is reversible and we can compare our (asymptotic) bound with the more precise one from \cite[Proposition 3]{Stroock-Diaconis}:
\[
\|P_i(n,\cdot) - \mu\|_{TV} \le \left(\frac{1-\pi(i)}{2\pi(i)}\right)^{1/2} \lambda_2^n 
= \frac1{\sqrt{2}} \times 0.3^n. 
\]
Clearly, the asymptotics of the bounds is the same for all three approaches.


\end{example}

\begin{example}[The new approach is similar to the MD one; the eigenvalue one is better]
Consider the MP with the state space $S = \{1,2,3\}$ and transition matrix
\[
{\cal P} = \left(
\begin{array}{l l l}
0 & 0.3 & 0.7  \\
0.3 & 0.7 & 0 \\
0.7 & 0 & 0.3
\end{array}
\right)\,.
\]
Here the matrix is doubly stochastic, so the distrbution $\pi=(1/3, 1/3, 1/3)$ is invariant; hence, clearly, the process is reversible: $\pi_ip_{ij} = \pi_jp_{ji}$ for all $i,j$; $\kappa(.) = 0.3$, $1-\kappa = 0.7$. 

Simple computations show that $V=0.7  * \hat {\cal P}$ with some stochastic matrix $\hat {\cal P}$; the latter  is a transition matrix for the process $\hat X_n$ (see the previous section). Hence, its spectral radius $r(V)$ equals 0.7. That is, our bound asymptotically coincides with the classical one. 

The characteristic equation is 
\[
\lambda^{3}-\lambda^{2}-0.37\lambda = 0,
\]
and so the eigenvalues of $\cal P$ (without hat) are $\lambda_{1}=1, \lambda_{2,3}=\pm \sqrt{37}/10 \approxeq \pm 0.6082763$. Hence, we have $|\lambda_{2,3}| < r(V)$.

\end{example}

\begin{example}[The new approach is similar to the classical (MD); the eigenvalue one is better]
\[
{\cal P} = \left(
\begin{array}{l l l}
0 & 0.3 & 0.7  \\
0.7 & 0 & 0.3 \\
0.3 & 0.7 & 0
\end{array}
\right) \, \mbox{-- is doubly stochastic;} \; \pi = (1/3, 1/3,1/3) \; \mbox{is invariant.}; \; \kappa(.)\equiv0.3.
\]
Therefore, $1-\kappa=0.7$. The operator $V$ has a form 
$$
V = 0.7 \times \hat {\cal P}, \quad \mbox{with some stochastic transition matrix $\hat {\cal P}$}. 
$$
It is known that such a matrix has its spectral radius $r(V)=0.7$. Thus, in this example our new bound is similar to the classical one. 
The analysis shows the following:
\[
|\lambda_{2,3}| = \sqrt{0.25 + 0.12} = \sqrt{0.37}  \approxeq 0.6082763 < 0.7.
\]
Note that the matrix $\cal P$ (without hat) is not reversible ($\pi_i p_{ij} = \pi_j p_{ji}$ fails: e.g., $p_{12}/3 \not = p_{21}/3$), so, formally, the more precise bound from \cite[Proposition 3]{Stroock-Diaconis} is not applicable. 
\end{example}

\begin{example}[The new approach is better than the MD and similar to the eigenvalue one]
Consider the MP with the state space $S = \{1,2,3\}$ and transition matrix
\[
{\cal P} = \left(
\begin{array}{l l l}
0 & 0.3 & 0.7  \\
1.0 & 0 & 0 \\
0.8 & 0.1 & 0.1
\end{array}
\right)\,.
\]
Here $\hat S^2 = \{(1,2), (1,3), (2,3), (2,1), (3,1), (3.2)\} =: (I,II, III, IV, V, VI)$, 
\[
\kappa(1,2) =\kappa(2,1) = 0,\;  \kappa(1,3) =\kappa(3,1) = 0.2,\; \kappa(2,3) =\kappa(3,2) = 0.8;
\]
\[
1-\kappa(1,2)= 1-\kappa(2,1) = 1,\; 
1-\kappa(1,3) =1-\kappa(3,1) = 0.8,\; 1-\kappa(2,3) =1-\kappa(3,2) = 0.2.
\]
The code -- to be published in the full version --  
provides the value $r(V)$: 
\[
1-\kappa=r(V) = 9/20+\sqrt{65}/20 \approx 0.85311289.
\]
The roots of the characteristic equation for the original transition probability matrix are $\lambda_{1}=1, \lambda_{2,3}=-9/20 \pm \sqrt{65}/20$, whilst the roots of the characteristic equation for the coupled process $\hat X_n$ are $\lambda_{c}=-9/20 - \sqrt{65}/20$, $\lambda_{c}=-9/20 + \sqrt{65}/20$, a repeated root $\lambda_{c}=0$, $\lambda_{c}=-\sqrt{65}/20 + 9/20$, $\lambda_{c}=\sqrt{65}/20 + 9/20$, where the last root has the highest value, i.e. $r=\sqrt{65}/20 + 9/20$ and so $|\lambda_{2}|=r=|-9/20 - \sqrt{65}/20| \approx 0.85311289$ indeed.
In this example the classical Markov -- Dobrushin method is useless since $1-\kappa=1$. The Remark \ref{rem1} is also applicable here.

\end{example}







\end{document}